\documentclass[12pt]{mymathart}
\usepackage{mymathmacros}

\numberwithin{equation}{section}

 \newtheoremstyle{numberedstyle}
   {9pt}
   {9pt}
   {\normalfont}
   {}
   {\bfseries}
   {.}
   {\newline}
   {}

\renewcommand{\dist}{\operatorname{dist}}

\newtheorem{thm}{Theorem}[section]%
\newtheorem{lem}[thm]{Lemma}%
\newtheorem{cor}[thm]{Corollary}%
\newtheorem{prop}[thm]{Proposition}%

\theoremstyle{numberedstyle}
\newtheorem{defn}[thm]{Definition}%

\renewcommand{\H}{\mathbb{H}}

\newcommand{\F}{\mathcal{F}}

\usepackage{hyperref}

\title[Exotic Baker and wandering domains]{Exotic Baker and wandering domains \\
    for Ahlfors islands maps}
\author{Lasse Rempe}
\address{Dept. of Math. Sciences, 
University of Liverpool, L69 7ZL, 
UK}
\email{l.rempe@liverpool.ac.uk}
\author{Philip J.\ Rippon}
\address{
   Department of Mathematics and Statistics \\
   The Open University \\
   Walton Hall\\
   Milton Keynes MK7 6AA\\
   UK}
\email{p.j.rippon@open.ac.uk}
\subjclass[2010]{30D05 (primary), 37F10, 30D45, 30E10 (secondary)}

\thanks{Both authors were supported by the European CODY network. The first author is supported by EPSRC fellowship EP/E052851/1. }

\begin{document}


 \begin{abstract}
 Let $X$ be a Riemann surface of genus at most $1$, i.e.\ $X$ is
  the Riemann sphere
  or a torus.
  We construct a variety of examples of analytic functions $g:W\to X$, where
   $W\subsetneq X$ is an arbitrary domain, that satisfy Epstein's
   ``Ahlfors islands condition''. In particular, we show that the accumulation
   set of any curve tending to the boundary of $W$ can be realized as the
   $\omega$-limit set of a Baker domain of such a function. As a corollary of our construction,
   we show that there are entire functions with Baker domains in which the iterates converge
   to infinity arbitrarily slowly. We also construct Ahlfors islands maps with wandering domains
   and logarithmic singularities, as well as examples where $X$ is
   a compact hyperbolic surface.
 \end{abstract}

 \maketitle

 \section{Introduction}

 The class of ``Ahlfors islands maps'' was introduced by
  Epstein, arising naturally from his concept of
  ``finite-type maps'' \cite{adamthesis}
  as a means of extending the iteration theory of
  transcendental entire or meromorphic functions to a more general
  setting. An Ahlfors islands map is an analytic function $g:W\to X$,
  where $W$ is
  an open subset of a compact Riemann surface $X$, that satisfies a certain
  maximality condition with respect to its domain of definition.
  This condition ensures that the
  Fatou set $F(g)$ and its complement, the Julia set
  $J(g)=X\setminus F(g)$, retain their usual properties.
  (See Section \ref{sec:definitions} for the formal definition. The theory of
   Ahlfors islands maps is developed in detail
   in the currently
   unpublished manuscripts \cite{finitetype,epsteinoudkerkahlfors}; compare
   also \cite{hypdim} for a short introduction.)

 If $g$ is an Ahlfors islands map, then---as
  for non-linear entire and meromorphic
  functions---any invariant component
  of the Fatou set $F(g)$ is one of finitely many types
  (see \cite[Chapter 5]{jackdynamicsthird} or \cite[Theorem~6]{waltermero}):
  an immediate basin of attraction for a (super-) attracting or parabolic
  fixed point,
  a rotation domain, or a \emph{Baker domain}; the latter is a periodic component of $F(g)$ in
  which the iterates converge to the boundary $\partial W$ of the domain
  of definition of $g$. If $U$ is a Baker domain or a \emph{wandering
  domain} (i.e., a component of $F(g)$ whose forward orbit is not eventually periodic),
  then by the \emph{$\omega$-limit set} $\omega(U)$ we mean
  the set of all points $w\in\cl{W}$ for which there is some point
  $z\in U$ whose orbit under $g$ accumulates on $w$.

 The $\omega$-limit set 
  of an invariant Baker domain is necessarily connected
  (see Lemma \ref{lem:accumulationsets}).
  Hence,
  if $X\setminus W$ is totally disconnected (as in the case of transcendental
  entire or meromorphic functions and their iterates) and $g$
  has a Baker domain $U$, then $\omega(U)$ consists
  of a single point.
  Otherwise, however, it is conceivable
  that a Baker domain of $g$ may have a nondegenerate continuum as its
  $\omega$-limit set. Epstein asked the question whether
  such ``exotic Baker domains'' can in fact exist.

 Oudkerk showed that the answer is ``yes'' 
  \cite{oudkerkexoticbakerdomains}.
  He constructed
  an example where $W$ is a simply-connected domain having a prime end with
  a nontrivial impression, and $g:W\to\Ch$ has a Baker domain
  whose $\omega$-limit  set is exactly this prime end impression.

 This raises a number of questions, e.g.\ whether the boundary of the
  domain of
  an Ahlfors islands
  map with an exotic Baker domain can be locally connected,
  or whether an Ahlfors islands map can have an exotic Baker domain whose
  $\omega$-limit set is $\partial W$. More generally,
  one can ask, given a domain $W$ and a compact connected set
  $K\subset\partial W$,
  whether
  there is an Ahlfors islands map with a Baker domain $U$ such
  that $\omega(U)=K$.

 In this note, we give a complete answer to this question
  in the case where $X$ is of genus $0$ or $1$. (That is,
  $X$ is the Riemann sphere, or a torus $X=\C/\Gamma$, where
  $\Gamma$ is a period lattice.)

 \begin{thm}[Exotic Baker domains] \label{thm:mainbaker}
  Let $X$ be a compact Riemann surface of genus at most $1$.
   Let $W\subsetneq X$ be nonempty,
   open and connected, and let $\gamma:[0,\infty)\to W$
   be an injective 
   curve with $\gamma(t)\to\partial W$ as $t\to\infty$. Let
   $K\subset \partial W$ be the accumulation set of~$\gamma$.

  Then there exists an Ahlfors islands map $g:W\to X$
   such that $g$ has a Baker domain $U$ with $\gamma\subset U$
   whose $\omega$-limit set is exactly $K$.

  If $X=\Ch$, then $g$
   can furthermore
   be chosen to omit any given value $a\in\Ch\setminus \gamma$.
 \end{thm}
 \begin{remark}[Remark 1]
  If an Ahlfors islands map $g$ has a Baker domain $U$, then
   $\omega(U)$ is the accumulation set of an injective
   $C^{\infty}$ 
   curve contained in $U$. (See Lemma \ref{lem:curvesinbakerdomains}.)
  Hence Theorem \ref{thm:mainbaker}
   states that any potential $\omega$-limit set is indeed
  realized by a Baker domain.
 \end{remark}
 \begin{remark}[Remark 2]
  If $X=\Ch$ and $W\subset\C$, then the final statement of the theorem
   implies that $g$ can be chosen to be analytic. In fact, this is
   exactly how the result is proved.
 \end{remark}

  Oudkerk's example is fairly explicit, and he
  verifies the Ahlfors islands property directly. Instead,
  we use Arakelian's approximation theorem (or, in the case of the torus,
  the generalization for Riemann surfaces
  due to Scheinberg \cite{scheinbergapproximation}) and
  the theory of normal families to obtain our examples. The main
  observation needed for the proof is that a suitable model function
  on a simply-connected domain tending to the boundary can be approximated
  sufficiently closely,
  even though by Nersesjan's theorem \cite[Satz IV.3.4]{gaier}
  \emph{arbitrarily}
  close
  approximation is not possible on such sets. See Lemma \ref{lem:approximation}
  for a precise statement.

 In \cite{ripponstallardslowescaping}
  the \emph{slow escaping set} of an entire (or meromorphic) function
  is studied.  In this context, it is interesting to note
  that we can control the growth of orbits inside the Baker domains we
  construct. In particular, we can construct entire functions
  having Baker domains in which the iterates tend to infinity
  arbitrarily slowly. Let us denote by $\dist_X$ distance with
  respect to the spherical
  metric (if $X=\Ch$) or with respect to the unique flat
  metric of area $1$ (if $X$ is a torus).

 \begin{prop}[Arbitrarily slow escape] \label{prop:speed}
  Let $\eps_n$ be an arbitrary sequence of positive numbers with
    $\lim_{n\to\infty} \eps_n = 0$. Then the function $g$ in Theorem
    \ref{thm:mainbaker} can be chosen such that the Baker domain
    $U$ also has the following property: for every $z\in U$, there
    exists $n_0\geq 0$ such that
       \[ \dist_X(g^n(z),\partial W)\geq \eps_{n} \]
    for all $n\geq n_0$.
 \end{prop}

 We also produce
  Ahlfors islands maps with wandering domains. As far as we know,
  they are the first such examples.

 \begin{thm}[Wandering domains] \label{thm:mainwandering}
  Let $X$ be a compact Riemann surface of genus at most $1$, and let
   $W\subsetneq \Ch$ be nonempty, open and connected. Also let
   $K\subset \partial W$ be compact.

  Then there exists an Ahlfors islands map $g:W\to X$
   such that $g$ has a wandering domain
   whose $\omega$-limit set is exactly $K$.

  If $X=\Ch$, then $g$ can furthermore be chosen to omit any given
   value $a\in\Ch$.
 \end{thm}
\smallskip


 Finally, we also produce examples of Ahlfors islands maps with
  logarithmic tracts.
  Recall that $a\in X$ is called a \emph{logarithmic asymptotic value}
  of a holomorphic
  function $g:U\to X$ if there is a simply-connected neighborhood $D$ of $a$
  and a component $T$ of $g^{-1}(D)$ such that $g:T\to D\setminus\{a\}$ is
  a universal covering. In this case, $T$ is called a \emph{logarithmic
  tract} of $g$.

 It was remarked in \cite{hypdim} that Arakelian's theorem can be used
  to construct Ahlfors functions in the unit disk having a logarithmic
  tract
  that spirals out to the unit circle.
  We show that in fact
   any injective $C^1$ curve tending to the boundary in a domain $W\subset X$
  (with $X$ as above) 
  can occur as an asymptotic curve for a logarithmic asymptotic value of
  an Ahlfors islands map $W\to X$.
  If this asymptotic value is itself in the accumulation set of the
  asymptotic curve, then by \cite[Theorem A.1]{hypdim}
  (which generalizes a result of \cite{baranskikarpinskazdunik})
  the Hausdorff
  dimension of the Julia set of $g$ is strictly greater than one.

 Since the existence of logarithmic asymptotic values is a function-theoretic
  rather than a dynamical property, it makes sense to state our
  results also for \emph{non-iterable} Ahlfors islands maps, where
  the domain $W$ is a nonempty open subset of some compact Riemann surface
  $Y$, which may or may not agree with the range $X$. 
  (See Remark 1 after Definition \ref{defn:ahlforsislands}.) 

  \begin{thm}[Logarithmic tracts] \label{thm:mainlogtracts}
   Let $X$ be a compact Riemann surface of genus at most~1 and let 
    $Y$ be a compact Riemann surface of arbitrary genus. Let
    $W\subsetneq Y$ be nonempty, open and connected.
    Also let $\gamma:[0,\infty)\to W$
    be an injective $C^1$ curve with $\gamma(t)\to\partial W$ as $t\to\infty$.

  Then for any $w_0\in X$ and any Jordan neighborhood $U\subset X$ of $w_0$,
   there exists an Ahlfors islands map
   $g:W\to X$ such that $g(\gamma)\subset U\setminus\{w_0\}$,
   $g(\gamma(t))\to w_0$ as
   $t\to\infty$, and such that the component of
   $g^{-1}(U\setminus\{w_0\})$ containing $\gamma$ is mapped as a universal covering by~$g$.

  If $X=\Ch$, then $g$ can furthermore be chosen to omit any given
   value $a\in\Ch$ with $a\notin\cl{U}\setminus\{w_0\}$.
 \end{thm}

 We note that, in particular, we obtain examples of Ahlfors islands
  functions with arbitrary domains whose range is a Riemann surface
  of genus at most $1$:

 \begin{cor}[Examples of Ahlfors islands maps]
 Let $X$ and $Y$ be compact Riemann surfaces, and let $W\subsetneq Y$ be
  nonempty, open and connected. If $X$ has genus at most $1$,
  then there exists an Ahlfors islands
  map $g:W\to X$. 
 \end{cor}
 If the surface $X$ is hyperbolic, then the question of which domains
  $W$ can support Ahlfors islands maps with values in $X$ is much more
  subtle. 
  For example, if $W\subset Y$ is simply-connected, then
  any universal covering $g:W\to X$ is 
  an Ahlfors islands map. On the other hand,
  an analytic map taking values in a hyperbolic surface cannot have
  isolated essential singularities. So if $\partial W$ has
  isolated points and $X$ is hyperbolic, then 
  there are no Ahlfors islands maps $g:W\to X$. 

 Hence we cannot expect to generalize our theorems to the case of
  arbitrary domains $W$ in a hyperbolic Riemann surface $X$. However,
  our methods can still be used to prove the following, slightly
  weaker results.

 \begin{thm}[Ahlfors maps on hyperbolic surfaces]\label{thm:mainhyperbolic}
   Theorems \ref{thm:mainbaker}, \ref{thm:mainwandering} and
    \ref{thm:mainlogtracts} still hold
    when $X$ is a compact hyperbolic Riemann surface, provided
    the phrase ``there exists an Ahlfors islands map $g:W\to X$'' is replaced by
    ``there exist a domain $W'\subsetneq W$ and an Ahlfors islands map
      $g:W'\to X$''. 
 \end{thm}

 As far as we know, these are the first known examples of Ahlfors islands
   maps with Baker domains or wandering domains
   on a surfaces of genus $g\geq 2$.

\subsection*{Acknowledgments}
  We thank Patricia Dom\'inguez, Adam Epstein and Gwyneth Stallard
  for interesting discussions.

 \section{Definitions and preliminaries} \label{sec:definitions}

 \subsection*{Basic notation}
   We denote the complex plane, the Riemann sphere, the
    unit disk and the right half plane by $\C$,
    $\Ch$, $\D$ and $\H$, respectively. The (Euclidean) disk of radius
    $\delta$ around $z\in\C$ is denoted by $B(z,\delta)$. We also
    denote Euclidean distance by $\dist_\C$. 

   Closures and boundaries will usually be taken in an underlying
    compact Riemann surface $X$ (which $X$ is meant should be
    clear from the context).
    Sometimes we consider relative closures, and it
    will be convenient to denote the relative closure of $A\subset W$ in
    $W\subset X$ by
        \[ \cl{A}_W := \cl{A}\cap W. \]

   We also denote hyperbolic distance on any
    hyperbolic Riemann surface $W$ by $\dist_W$. Then
    $B_W(z,\delta)$ denotes the hyperbolic disk in $W$ of radius $\delta$
    around $z$. (See e.g.\
    \cite[Chapter 2]{jackdynamicsthird} for basic definitions and
    results of hyperbolic geometry, and in particular \emph{Pick's theorem},
    which we shall use frequently.)

 \subsection*{Ahlfors islands maps and normal families}
  For the following definitions, let $X$ and $Y$
   be compact Riemann surfaces and
   let $W\subset Y$ be open and nonempty.
   It will be convenient to introduce the following terminology.
   If $V\subset Y$ is a connected open set that intersects the boundary
   of $W$ in $Y$, then we shall call a component $U$ of
   $V\cap W$ a \emph{boundary neighborhood} of $W$, provided that $U$ is
   hyperbolic. The latter will always be the case if $V$ is chosen
   sufficiently small to omit at least three points of $Y$.

  A \emph{Jordan domain} on $X$ is a simply-connected domain
   that is bounded by a Jordan curve in $X$. If $V\subset X$ is a Jordan domain
   and $g:W\to X$ is holomorphic,
   then a \emph{simple island} over $V$ is a domain $I\subset W$ such that
   $g:I\to V$ is a conformal isomorphism.

 \begin{defn}[Ahlfors islands maps] \label{defn:ahlforsislands}
  A holomorphic function $g:W\to X$ has the
   \emph{Ahlfors islands property} and $g$ is called an \emph{Ahlfors islands map} if there exists a number $k$ such that
   the following is true.

   If $V_1,\dots,V_k$ are Jordan domains in $X$ with pairwise disjoint
    closures, then for every boundary neighborhood $U$ of $W$, there exists
    $j\in\{1,\dots,k\}$ such that $U$ contains a simple island
    of $g$ over $V_j$.
 \end{defn}
 \begin{remark}[Remark 1] 
   In \cite{hypdim}, the definition was stated only for the case 
    $X=Y$, in which case  it is possible to consider $g$ as a
    dynamical system. However, this is not
    required and the above is the general definition as introduced by Epstein.
    For our purposes,
    the additional flexibility simplifies the discussion when
    constructing Ahlfors islands maps
    on the torus.
 \end{remark}
 \begin{remark}[Remark 2]
  One of the key properties of the definition is that the composition of
   two Ahlfors islands maps is again an Ahlfors islands map. We do not
   require this fact.
 \end{remark}

 If $X=Y$ and $g:W\to X$ is an Ahlfors islands map,
  then the \emph{Fatou set}
  consists of all points $z\in X$ for which there exists either
  some $k\geq 0$ such that $g^k(z)\in X\setminus\cl{W}$ or an open neighborhood
  $U$ of $z$ in which the iterates $g^k|_U$ are all defined and form a
  normal family. By virtue of the Ahlfors islands property of $g$,
  the Julia set $J(f) = X\setminus F(f)$ is the closure of the set of
  repelling periodic points, and also retains its other well-known
  properties. In particular, if $g$ is \emph{non-elementary}, i.e.\
  not a conformal isomorphism of $X$, then $J(g)$ is a nonempty, perfect,
  compact set. (We do not use these facts in this article.)

 \smallskip

 The definition of Ahlfors islands maps---and their name---%
  is inspired by the classical five islands theorem
  of Ahlfors, which implies that every transcendental meromorphic function
  is an Ahlfors islands map (with $k=5$). We will use the following
  ``normal families
  version'' of this result; see e.g.\ \cite{walterahlfors}.

 \begin{thm}[Ahlfors islands theorem, normal families version]
    \label{thm:ahlforsnormal}
  Let $V_1,\dots,V_5\subset\Ch$ be Jordan domains with pairwise disjoint
   closures. If
   $U\subset\C$ is a domain, and $\F$ is a family of
   meromorphic functions $f:U\to \Ch$ that have no simple islands over any
   of the $V_j$, then $\F$ is a normal family. \qedd
 \end{thm}
 \begin{remark}
  If the functions in $\F$ are holomorphic---which is the case
   we are interested in---the number five can be replaced by four.
   (An extremal example is given by the iterates of the sine function,
    which has no islands over any Jordan domain that contains $1$,
    $-1$ or $\infty$.)
 \end{remark}

 We will use Theorem \ref{thm:ahlforsnormal} to construct our examples
  by introducing a class of functions that are not normal near any
  point of the boundary of $W$.
  Note that such functions can exist only if $X$ is an elliptic or parabolic
  manifold (i.e., $X$ is the sphere or a torus), since any family of
  analytic functions taking values in a hyperbolic
  surface is normal.

 \begin{defn}[Strong non-normality] \label{defn:strongnonnormal}
  Suppose that $W$ is a hyperbolic Riemann surface.

   We shall say that
    a holomorphic function $g:W\to \C$ is \emph{strongly non-normal}
    if there exists a number $\delta>0$ with the following property:

  Whenever $w_n\in W$ is a sequence tending to the boundary of $W$ and
   $\pi_n:\D\to W$ are universal covers with $\pi_n(0)=w_n$, the
   family
    \[ \{ g\circ \pi_{n}: B_{\D}(0,\delta)\to\Ch \} \]
   is not normal.
 \end{defn}

 Suppose that $g:W\to \C$ is strongly non-normal, where
  $W$ is a hyperbolic subdomain of a compact Riemann surface $Y$.
  Since
  every boundary neighborhood of $W$
  contains a sequence of hyperbolic balls of
  fixed diameter tending to the boundary, it follows from the
  normal families version of the Ahlfors islands theorem that
  $g:W\to\Ch$ is an Ahlfors islands map.

 As our definition of strong non-normality
  is phrased in conformally invariant terms, it is
  invariant under precomposition with a conformal map; this is why we have
  chosen it among a number of other definitions that would also
  be suitable for our purposes. We will use the following condition
  to ensure that a function is strongly non-normal.

 \begin{lem}[Sufficient condition for strong non-normality]
     \label{lem:stronglynormal}
  Let $W$ be a hyperbolic Riemann surface, and let
   $g:W\to\C$ be holomorphic. Suppose that there are relatively
   closed subsets $B,C\subset W$ and a positive constant $\Theta$
   with the following properties:
   \begin{itemize}
     \item for every $w\in W$, $\dist_{W}(w,B)\leq \Theta$ and
         $\dist_{W}(w,C)\leq \Theta$;
     \item $g|_B$ is bounded;
     \item $g(c)\to\infty$ as $c\to\partial W$ within $C$.
   \end{itemize}

  Then $g$ is strongly non-normal.
 \end{lem}
 \begin{proof}
  Set $\delta := 2\Theta$; then every closed
   hyperbolic disk of radius $\delta/2$
   intersects both $B$ and $C$.
   Let $w_n$ and $\pi_n$ be as in the definition of strong non-normality.
   Then there are sequences $b_n$ and $c_n$ in $\cl{B_{\D}(0,\delta/2)}$
   such that $\pi_{n}(b_n)\in B$ and $\pi_{n}(c_n)\in C$,
   and in particular
         \[ \limsup_{n\to\infty} |g(\pi_{n}(b_n))|<\infty, \quad g(\pi_{n}(c_n))\to \infty\text{ as } n\to\infty. \]
   Clearly the sequence $g\circ \pi_{n}$ cannot be normal on
   $B_{\D}(0,\delta)$, as required.
 \end{proof}

 \subsection*{The Arakelian-Scheinberg approximation theorem}

 We use the following terminology from \cite{gaier}.
 \begin{defn}[Weierstra{\ss} sets]
  Let $Y$ be a compact Riemann surface, and let
   $W\subsetneq Y$ be a nonempty open and connected set.
   A relatively closed set $A\subset W$ is called a
   \emph{Weierstra{\ss} set} (in $W$)
   if any continuous function
   $f:A\to\C$ that is holomorphic on the interior of $A$ can be
   uniformly
   approximated,
   up to an arbitrarily small error $\eps$, by a holomorphic function
   $g:W \to\C$.
 \end{defn}

 \begin{thm}[Arakelian-Scheinberg approximation
   theorem] \label{thm:arakelian}
   A relatively closed set~$A$ in $W$ is a Weierstra{\ss} set
   if and only if $\hat{W}\setminus A$ is connected and locally connected,
    where $\hat{W}=W\cup\{\infty\}$ is the one-point compactification of $W$.
 \end{thm}
 \begin{remark}
  The case where $Y=\Ch$ is \emph{Arakelian's theorem}
  \cite[Satz IV.2.3]{gaier}; the general case is due to
  Scheinberg \cite{scheinbergapproximation}.
  The latter more generally treats \emph{arbitrary} non-compact
   Riemann surfaces $W$. In this setting, it turns out that there is no
   longer a topological characterization of Weierstra{\ss} sets, but
   Scheinberg gives a sufficient condition on the topology of $A$ under which
   the above condition is still necessary and sufficient. This condition
   includes the case where $W$ has finite genus.
 \end{remark}

  We note that it follows from Theorem \ref{thm:arakelian} that
   any countable disjoint union of nonseparating compact subsets
   $K_n\subset W$ tending to $\partial W$ as $n\to\infty$ is a
   Weierstra{\ss} set.

 Our construction will rely on approximating certain functions
  up to an error prescribed by
  a given function. Nersesjan's theorem \cite[Satz IV.3.4]{gaier} gives
  a necessary and sufficient condition
  on the set $A$ such that this is
  possible with an arbitrary error function in the case where
  $W\subset\Ch$. Unfortunately
  this criterion excludes, in particular, the case where $A$ has unbounded
  interior components.

 Instead, we will make use of the following trick
    \cite[Hilfssatz IV.3.3]{gaier} that
  also goes back to Arakelian.

 \begin{lem}[Approximation up to an error function]
     \label{lem:arakelian}
  Let $Y$ be a compact Riemann surface, let
   $W\subsetneq Y$ be nonempty, open and connected, and let
   $A$ be a Weierstra{\ss} set in $W$.
   Furthermore, let $f:A\to \C$ and
   $E:A\to\C$ be continuous on $A$ and holomorphic on the interior of $A$.

   Then there is a holomorphic function $g:W\to\C$ such that
         \[ |g(z) - f(z)| \leq |\exp(E(z))| \]
    for all $z\in A$.
 \end{lem}
 \begin{proof}
   Because $A$ is a Weierstra{\ss} set, we can find a holomorphic function
     $H:W\to\C$ such that
      \[ |H(z) - (E(z) - 1)| \leq 1 \]
    for all $z\in A$. In particular, we have $\re H(z)\leq \re E(z)$, and hence
      \[ |\exp(H(z))| = \exp(\re H(z)) \leq
            \exp(\re E(z)) = |\exp(E(z))| \]
    for all $z\in A$.

    Now let $G:W\to\C$ be a holomorphic function that approximates
     $f(z)\cdot \exp(-H(z))$ up to an error of at most $1$ on $A$. Then
     the map
     \[ g(z) := G(z)\cdot \exp(H(z)) \]
     satisfies
      \begin{align*}
          |g(z) - f(z)| &= |\exp(H(z))|\cdot |G(z) - f(z)\cdot \exp(-H(z))|
            \\ &\leq
         |\exp(H(z))| \leq |\exp(E(z))|
      \end{align*}
      for all $z\in A$.
 \end{proof}
 \begin{remark}
  Note that Lemma \ref{lem:arakelian} implies that approximation up to
   an arbitrary error function is possible if $A$ is a Weierstra{\ss} set that
   has only compact connected
   components; this is a special case of the
   theorem of Nersesjan.
 \end{remark}

 In Lemma \ref{lem:approximation}, we will use Lemma
  \ref{lem:arakelian} to approximate
  a function given inside a simply-connected subset of $W$ that accumulates
  on $\partial W$.

\subsection*{Limit sets}

 For completeness, we note two simple facts about $\omega$-limit sets.
 \begin{lem}[Limit sets] \label{lem:accumulationsets}
  Let $X$ be a compact Riemann surface, let $W\subset X$ be open and nonempty,
   and let
   $g:W\to X$ be an Ahlfors islands map.

  Also let $U$ be a Baker domain or a wandering domain of $g$. If $z\in U$
   and $(n_k)$ is a sequence such that $g^{n_k}(z)\to w_0\in\omega(U)$, then
   $g^{n_k}|_U\to w_0$ locally uniformly.

  In particular, $\omega(U)$ coincides with the set of limit points
   of the sequence
   $(g^n(z))$ for every $z\in U$.

  If $U$ is an invariant Baker domain, then there is a curve 
   $\alpha\subset U$ whose accumulation
   set is $\omega(U)$. In particular, $\omega(U)$ is connected.
 \end{lem}
 \begin{proof}
  Set
   $\wt{U} := \bigcup_{j\geq 0}g^n(U)$. Then, by assumption on $U$, every
   limit function of the normal family $(g^n|_U)$
   takes values in the nowhere dense set
   $\partial\wt{U}$, and hence is constant. In particular, with notation
   as in the lemma, every
   limit function of the sequence $(g^{n_k}|_U)$ is the constant function
   $w_0$. This proves the first claim.

  If $U$ is an invariant Baker domain and $z\in U$ is arbitrary, then we can pick
   some curve $\alpha_0\subset U$ connecting $z$ and $g(z)$. Then
   $\alpha:=\bigcup_{j\geq 0}g^j(\alpha_0)$ is the desired curve for the final claim.
 \end{proof}

 We can make the final part of the preceding lemma more precise as follows.

 \begin{lem}[Curves in Baker domains] \label{lem:curvesinbakerdomains}
   Let $X$ be a compact Riemann surface, and let $U\subset X$ be a hyperbolic
   domain.

   Let $f:U\to U$ be holomorphic, and suppose that
    $\dist_X(f^n(z),\partial U)\to 0$
    as $n\to \infty$. Let $z\in U$ and let
    $K$ be the set of limit points of the sequence
    $(f^n(z))$.

   Then there exists an injective $C^{\infty}$ curve $\gamma:[0,\infty)\to U$
    whose accumulation set is $K$.
 \end{lem}
 \begin{remark}
  We state this lemma only for completeness; it will not be used in
   the rest of the article.
 \end{remark}
 \begin{proof}
  Let
   $C$ be the hyperbolic distance between $z$ and $f(z)$. By
   Pick's theorem, $C$ is an upper bound for the hyperbolic distance
   between $f^n(z)$ and $f^{n+1}(z)$, for all $n$. So
     \[ V := \bigcup_{n\geq 0} B_U(f^n(z),2C) \]
   is a connected subdomain of $U$. Since the diameter of
    $B_U(f^n(z),2C)$ in $X$ tends to zero as $n\to\infty$, we have
    $\partial V\cap \partial U = K$.

  Now pick an  injective
   $C^{\infty}$ curve $\alpha:[0,\infty)\to V$ that
  tends to $\partial U$ as $t\to\infty$ and passes through all
  $f^n(z)$. Then the accumulation set of $\alpha$ is precisely $K$ (it
  contains $K$ since it passes through all $f^n(z)$
  and is contained in $K$ by choice of $V$).

  (That the existence of such a curve $\alpha$ follows from standard methods
   of topology should be plausible to the reader,
   but we sketch a proof for completeness.
   First note that it is sufficient to prove this statement
   without the requirement that $\alpha$ be $C^{\infty}$,
   as we can approximate any
   injective curve by a smooth and injective one.
   As in the previous lemma, there is a curve $\beta$ with the desired
   properties, but $\beta$ is not necessarily injective. This curve has an
   injective subcurve $\beta_0$ that tends to $\partial U$, but might not
   contain all the points
   $f^n(z)$. This is easily corrected by connecting the missing points,
   one by one, to $\beta_0$ using a piece of the curve $\beta$, modifying the
   result so as to obtain an injective curve.)
  \end{proof}

 \section{Wandering domains} \label{sec:wandering}

  \begin{proof}[Proof of Theorem \ref{thm:mainwandering}]
   There are three different cases to consider:
   \begin{enumerate}
    \item $X=\Ch$ and $W$ is a plane or punctured plane.
       In this case, we may assume without loss of generality that
       $W=\C$ or $W=\C\setminus\{0\}$. \label{item:parabolic}
    \item $X=\Ch$ and $W$ is a hyperbolic domain.
    \item $X=\C/\Gamma$ is a torus and $W$ is a hyperbolic domain.
   \end{enumerate}

   The proof is similar for all three cases, so we shall treat them in
    parallel, and remark on differences in the appropriate places.
    (We note that case (\ref{item:parabolic}) is well-known, but since many
     of the arguments will reoccur in later proofs, we include a proof for
     the parabolic case here for completeness.)

   If $X=\Ch$, then let
    $\pi:\Ch\to\Ch$ be a M\"obius transformation with
    $\pi(\infty)=a$ (where $a$ is the value that our map $g$ will omit; recall the statement of
    the theorem). Otherwise, $X=\C/\Gamma$ is a torus, and we let $\pi:\C\to X$ be the natural
    projection.

   The basic structure of the proof is as follows. We pick a
    sequence $z_k\in W$, $k\in\N$, that accumulates exactly on the given compact set
    $K\subset\partial W$ (and tends to the boundary ``sufficiently quickly'').
    We also pick a sequence of small Jordan domains $\Delta_k$
    around $z_k$, such that $\cl{\Delta_k}\subset W$ and the $\cl{\Delta_k}$ are pairwise disjoint,
    and let $\Delta_k'\subset\C$ be a component of
    $\pi^{-1}(\Delta_k)$. If $X=\Ch$, we assume that $z_k$ and $\Delta_k$
    are chosen so that $a\notin \cl{\Delta_k}$ for all $k$. Thus every
    $\Delta_k'$ is a bounded Jordan domain in $\C$.

   Then we choose a function $f:\bigcup_{k=1}^{\infty}\Delta_k\to\C$ with
    $f(\cl{\Delta_k})\subset \Delta_{k+1}'$ for all $k$, and approximate this
    map, using Lemma
    \ref{lem:arakelian}, by an analytic function $h:W\to\C$
    that still takes $\cl{\Delta_k}$ into
    $\Delta_{k+1}'$. This construction is done in such a way
    that we ensure that $h$ (considered as a function $W\to\Ch$) is an
    Ahlfors islands map.

   Because $\pi$ is a covering map onto $X$,
    the function $g := \pi\circ h:W\to X$ is clearly an
    Ahlfors islands
    map also. We have $g(\cl{\Delta_k})\subset \Delta_{k+1}$, and if $z_k$ and
    $\Delta_k$ are suitably chosen, then we can ensure that each $\Delta_k$
    is contained in a wandering domain of $g$, completing the proof.
    Note that, if $X=\Ch$, then $g$ omits the value $a$.

  We now provide the details.

 \smallskip

   \emph{Choice of $z_k$ and $\Delta_k$ in the hyperbolic case.}
    First suppose that $W$ is a hyperbolic domain.
     Then we may suppose that the sequence $z_k$ mentioned above
     is chosen such that
     \begin{equation}
       \dist_{W}(z_{k+1} , \{ z_1, \dots, z_k \}) \to \infty\;\text{ as } k\to\infty.
            \label{eqn:hypdist}
     \end{equation}
    Our only additional requirement on our sequence $\Delta_k$
     of Jordan domains in
     this case is that the hyperbolic diameter of $\Delta_k$ should
     be bounded by some constant $M$, independently of $k$.

 \smallskip

    \emph{Choice of $z_k$ and $\Delta_k$ in the parabolic case.}
     If $W=\C$ or $W=\C\setminus\{0\}$, we let $\bigl|\log|z_k|\bigr|$
      grow sufficiently
      quickly and pick $\Delta_k$ sufficiently small so that
        \[
           \bigl|\log|z|\bigr| > k\cdot \bigl|\log|w|\bigr|
            \]
     whenever $z\in \Delta_{k+1}$ and $w\in \Delta_k$.

 \smallskip

    \emph{Choice of $B$ and $C$.}
     In addition, we choose two discrete countable sets $B,C\subset W$,
      disjoint from each other and all $\cl{\Delta_k}$,
      such that the accumulation sets of $B$ and of $C$ in $X$ coincide
      exactly with $\partial W$. (These sets will be used to ensure that
      the function $h$ is an Ahlfors islands map.)

    If $W$ is hyperbolic, then we also require that these sets are chosen
     such that points of
     $W$ have uniformly
     bounded hyperbolic distance from both $B$ and $C$ (as in
     Lemma \ref{lem:stronglynormal}). Note that this is possible
     because we required the $\Delta_k$ to have uniformly bounded
     hyperbolic diameters.

 \smallskip

  \emph{Definition of $f$ and approximation.}
    We set
       \[ A := B\cup C \cup \bigcup_{k=1}^{\infty} \cl{\Delta_k}. \]
    Observe that $A$ is a Weierstra{\ss} set by Theorem \ref{thm:arakelian}.

    We now define a continuous function
     $f:A\to\C$ as follows:
     \begin{enumerate}
      \item $f|_{\Delta_k}$ is a holomorphic function that extends continuously
        to $\partial \Delta_k$ such that $f(\cl{\Delta_k})\subset \Delta_{k+1}'$.
        (Recall that $\Delta_{k+1}'$ is a component of $\pi^{-1}(\Delta_{k+1})$.)
            \label{item:fonwanderingdomain}
      \item
       $f|_B$ is bounded.  \label{item:bconstant}
      \item $f(c)\to \infty$ as $c\to\partial W$ within $C$.
      \end{enumerate}

    We also define a locally constant function $e(z):A\to(0,\infty)$ by setting
     $e\equiv 1$ on $B$ and $C$, and
       \[ e(z) := \dist_{\C}(f(\Delta_k),\partial \Delta_{k+1}') \]
     on $\cl{\Delta_k}$.

    Now we can apply Lemma \ref{lem:arakelian} to $f$ and
     $E = \log e$ to obtain a holomorphic function $h:W\to\C$ such that
       \[ |h(z) - f(z)| \leq e(z) \quad \text{for all } z\in A. \]

   Set $g:= \pi\circ h$. We have $g(\Delta_{k})\subset
      \pi(\Delta_{k+1}') = \Delta_{k+1}$, so $\Delta_k$ is contained
     in the Fatou set of $g$ for all $k$. Furthermore, the iterates of
     $g$, restricted to $\Delta_k$, tend to $\partial W$, so
     the Fatou component containing $\Delta_k$ is either a
     a Baker domain (or pre-image component of a Baker domain) or
     a wandering domain.

 \smallskip

   \emph{$h$ is an Ahlfors islands map.} If $W$ is
     hyperbolic, then it follows from Lemma \ref{lem:stronglynormal} and
     our definitions that $h$ is strongly normal, and hence an
     Ahlfors islands map.

    If $W$ is parabolic, then it follows from the fact that
     $h$ is bounded on $B$ and tends to infinity on $C$ that
     each point of $\partial W$ is
     an essential singularity of $h$. Hence $h$ is an Ahlfors islands map
     by the classical five islands theorem.

    As noted above, $g$ is then also an Ahlfors islands map.

  \smallskip

   \emph{Each $\Delta_k$ is contained in a wandering domain.}
    First suppose that $W$ is hyperbolic. If $\Delta_1$ was contained in
     an eventually periodic domain, then
     there would be some
     $k\in\N$ and $l>k$ such that $z_k$ and $z_l$ can be connected by
     a curve in the Fatou set $F(g)$. Recall that the
    diameter of $\Delta_k$ in the hyperbolic distance of
     $W$ is at most $M$. Hence we would have,
     by Pick's theorem,
     that
      \begin{align*}
         \dist_{F(f)}(z_k,z_l) &\geq
         \dist_{F(f)}(g^m(z_k),g^m(z_l)) \geq
         \dist_{W}(g^m(z_k),g^m(z_l)) \\ &\geq
         \dist_{W}(z_{k+m},z_{l+m}) - 2M \to \infty \end{align*}
    as $m\to \infty$ by (\ref{eqn:hypdist}). This is a
     contradiction.

    If $W$ is parabolic, then we can likewise show that $\Delta_k$ must be
     contained in a wandering domain. Indeed, since
     $g(\Delta_k)\subset \Delta_{k+1}$, we have
      \begin{equation}
         \label{eqn:wanderinggrowth}
          \bigl|\log |g^{k+1}(z)|\bigr| \geq k\cdot \bigl|\log|g^k(z)|\bigr|
      \end{equation}
     for all $z\in \Delta_1$.

   It is well-known that there is a bound
     on the speed of escape in Baker domains \cite[Lemma 7]{waltermero}.
     More precisely, if $z$ is eventually mapped to a Baker domain
     of period $l$, then for
     sufficiently large $k$ and some suitable constant $C$, we have
         \[ \bigl|\log |g^{k+l}(z)|\bigr| \leq C\cdot \bigl|\log|g^k(z)|\bigr|.
                      \]
     But this would contradict (\ref{eqn:wanderinggrowth}).
      Thus $\Delta_1$ is contained in a wandering domain, and hence so is
      each $\Delta_k$.

  By construction and Lemma
   \ref{lem:accumulationsets}, 
   the $\omega$-limit set of these wandering domains is exactly
     $K$, as claimed.
    \end{proof}

 \section{Baker domains}  \label{sec:baker}

  To prove Theorem \ref{thm:mainbaker}, we
   begin by forming a simply-connected domain around a given
   curve $\gamma$. The following lemma shows that this is always possible.
   This result is surely known, but we do not know of a reference, and
   therefore include a proof in Section~\ref{sec:strip}.

 \begin{lem}[Simply-connected domain around a curve] \label{lem:strip}
  Let $W$ be a noncompact
   Riemann surface, with a metric $d$ on $W$ (compatible with
   the topology). Furthermore, let $\gamma:[0,\infty)\to W$ be an
   injective curve such that $\gamma(t)\to\infty$ as
   $t\to\infty$ (in the one-point compactification $\hat{W}$).
   Also let $\delta:[0,\infty)\to (0,\infty)$ be continuous.

  Then there exists a simply-connected domain
   $V\subset W$ with $\gamma\subset V$ and
   a conformal isomorphism $\phi:\H\to V$ such that
   \begin{enumerate}
    \item $\displaystyle{V\subset\bigcup_{t\geq 0}%
                  \{z\in W: d(z,\gamma(t))<\delta(t)\}}$.
    \item $\phi$ extends continuously to a homeomorphism between the closures
      $\cl{\H}_{\C}$ and $\cl{V}_W$ of $\H$ and $V$ in $\C$ resp.\ $W$.
    \item $\phi(z)\to\infty$ in $\hat{W}$ as $z\to\infty$ in $\H$.
    \item If $\alpha\subset V$ is any curve that tends to $\infty$ in
      $\hat{W}$, then points of $\alpha$ have uniformly bounded distance
      to $\gamma$, and vice versa.  \label{item:distancebetweencurves}
  \end{enumerate}
 
  If $\gamma$ is $C^1$, then furthermore the domain $V$ 
   can be chosen such that
  \begin{enumerate}
    \item[(e)] $\re\phi^{-1}(\gamma(t))\to \infty$ as $t\to\infty$.
         \qedd
   \end{enumerate}
 \end{lem}
 \begin{remark}[Remark 1]
  We do not use (e) in the proof of Theorem \ref{thm:mainbaker}, but
   will require it in the next section to prove
   Theorem~\ref{thm:mainlogtracts}.
 \end{remark}
 \begin{remark}[Remark 2]
  Suppose that (under the hypotheses of Lemma \ref{lem:strip}),
   $W$ is a hyperbolic subdomain of
   some compact Riemann surface $Y$,
   $d$ is the hyperbolic metric on $W$ and $\delta$ is bounded from above.
   Then (\ref{item:distancebetweencurves}) implies that
    any curve $\alpha\subset V$ tending to $\partial W$ has the
    same accumulation set as $\gamma$.
 \end{remark}

 Given a domain $V$ as in the lemma, we then aim to
   construct a holomorphic function $g:W\to X$ that maps $V$ into itself.
   To do so, we require
   a version of the approximation theorem
   that allows us to control the error in the (potentially very thin) strip
   $V$. The following result provides the kind of control that we are looking
   for.

  \begin{lem}[Control of approximation on simply-connected sets]
    \label{lem:approximation}
   Let $X$ be a compact Riemann surface, let
    $W\subsetneq X$ be a domain, and let $A\subset W$ be a Weierstra{\ss} set.
    Suppose that $V\subset \C$ is a simply-connected domain, let
    $\phi:\H\to V$ be a Riemann map, and set 
    $V' := \phi(\{\zeta\in\H: \re \zeta \geq 1\})$.

   Now let $f:A\to\C$ be continuous on $A$ and holomorphic in the interior
    of $A$. Suppose furthermore that $A$ can be written as the disjoint
    union of two relatively closed subsets $A_1$ and $A_2$, and that
    $f(A_1)\subset V'$.

 Then for every
    $\eps>0$, there exists a holomorphic function
    $g:W\to\C$ with the following properties.
    \begin{enumerate}
      \item $|f(z) - g(z)|\leq \eps$ for $z\in A_2$.
      \item For all $z\in A_1$,
       $g(z)\in V$ and $\displaystyle{%
       |\phi^{-1}(f(z)) - \phi^{-1}(g(z))| \leq \eps}$.\label{item:approxinV}
    \end{enumerate}
  \end{lem}
  \begin{proof}
   Consider the function
     \[ \Psi : V \to \C\setminus\{0\};
       z \mapsto \lambda\cdot\left( z - \phi\left(\phi^{-1}(z)+\frac{1}{2}\right)\right), \]
    where $\lambda>0$ is a constant. We claim that $\lambda$ can be chosen
    such that $B(z,|\Psi(z)|)\subset V$ for all $z\in V'$ and furthermore
     \begin{equation}
        \phi^{-1}( B(z,|\Psi(z)|) ) \subset
           B(\phi^{-1}(z),\eps). \label{eqn:containment}
     \end{equation}
    Indeed, let $z\in V'$ and set $\zeta := \phi^{-1}(z)$.
    By K\"obe's theorem, there is a constant
    $C$ such that
      \[ \left|\phi(\zeta) - \phi\bigl(\zeta + \frac{1}{2}\bigr)\right| \leq
              C |\phi'(\zeta) | \]
     and furthermore
      \[ \phi(B(\zeta,\eps)) \supset B\bigl(\phi(\zeta),\min(1,\eps)\cdot |\phi'(\zeta)|/4\bigr). \]
     Hence, if $\lambda < \min(1,\eps)/4C$, the claim follows.

    Since $V$ is simply-connected and $\Psi(z)\neq 0$, there exists
     a holomorphic logarithm
      $\psi(z): V\to \C$ of $\Psi$; i.e., $\exp(\psi(z)) = \Psi(z)$ for $z\in V$.

    We define
       \[ E:A\to\C; \quad z\mapsto\begin{cases}
                          \psi(f(z)) & z\in A_1 \\
                          \log\eps & z\in A_2. \end{cases} \]
     Since $A_1$ and $A_2$ are relatively closed, the function
     $E$ is continuous on $A$ and holomorphic in the interior of $A$.
     Applying Lemma \ref{lem:arakelian} to $f$ and $E$, we obtain a
     holomorphic function $g:W\to \C$ such that
        \[ |g(z) - f(z)| \leq |\exp(E(z))| \quad\text{for } z\in A. \]

    So we have $|g(z)-f(z)|\leq \eps$ on $A_2$. Furthermore, suppose that
     $z\in A_1$, so $f(z)\in V'$ by hypothesis. Then, by the above,
      $g(z)\in B(f(z),|\Psi(f(z))|)\subset V$ and
      \[ \phi^{-1}(g(z)) \in \phi^{-1}(B(f(z),|\Psi(f(z))|))
             \subset B(\phi^{-1}(f(z)),\eps) \]
     by (\ref{eqn:containment}). This completes the proof.
  \end{proof}

 We are now ready to prove Theorem \ref{thm:mainbaker} and
   Proposition \ref{prop:speed}.

  \begin{proof}[Proof of Theorem \ref{thm:mainbaker}]
   Let $\pi$ be defined as in the proof of Theorem \ref{thm:mainwandering};
    i.e. $\pi:\C\to X$ is a projection when $X$ is a torus; otherwise
    $\pi:\Ch\to X$ is a M\"obius transformation taking $\infty$ to
    a given point $a\in \Ch\setminus\gamma$.
   Also pick discrete sets $B$ and $C$, disjoint from $\gamma$,
    as in the proof of Theorem
    \ref{thm:mainwandering}. 

   We now apply Lemma \ref{lem:strip}, taking $d$ to be 
    a flat metric (if $W$ is the plane or the punctured
    plane) or  
    the hyperbolic
    metric in $W$ (otherwise). 
    The function $\delta$ is any continuous function into the positive
    reals with the property that
    $d(\gamma(t),B\cup C)> \delta(t)$ for all $t\in [0,\infty)$.

   Let $V$ be the simply-connected domain given by Lemma \ref{lem:strip},
    and let
    $\phi:\H\to V$ be the corresponding conformal isomorphism, whose
    continuous extension to $\cl{\H}_{\C}$ we also denote by
    $\phi$. Let $\tilde{V}$ be a component of $\pi^{-1}(V)$; then
    $\tilde{V}\subset\C$ is a simply-connected domain and
    $\pi:\tilde{V}\to V$ is a conformal isomorphism.
    Let $\tilde{\phi}: \H\to\tilde{V}$ be the conformal isomorphism
    satisfying $\pi\circ\tilde{\phi} = \phi$.

  We will use
    Lemma \ref{lem:approximation} to construct a strongly non-normal function
    $h:W\to\C$ that maps $\cl{V}_W$ into $\tilde{V}$.

    To do so, set $A := \cl{V}\cup B \cup C$. Note that $A$ is
     a Weierstra{\ss} set by Theorem \ref{thm:arakelian}. Define
     a function $f:A\to \C$, continuous on $A$ and holomorphic in
     the interior of $A$, such that
    \begin{enumerate}
      \item $f(z) := \tilde{\phi}(\phi^{-1}(z) + 2 )$ for $z\in \cl{V}$;
      \item $f$ is bounded on $B$;
      \item $f(c)\to \infty$ as $c\to\partial W$ within $C$.
    \end{enumerate}

    We also set $A_1 := \cl{V}$ and $A_2 := B\cup C$. Then we can
     apply Lemma \ref{lem:approximation} (with some $\eps<1$), to
     obtain a holomorphic function $h:W\to \C$. As in the proof of
     Theorem \ref{thm:mainwandering}, $h$ is an Ahlfors islands map,
     as is $g := \pi \circ h$.

    We have $h(\cl{V}_W)\subset \tilde{V}$ by Lemma
     \ref{lem:approximation} (\ref{item:approxinV}). Let us define
     $G:\cl{\H}_{\C}\to\H$ by
     \[ G(\zeta) := \phi^{-1}(g(\phi(\zeta))) = \tilde{\phi}^{-1}(h(\phi(\zeta))). \]
    Then, for any $\zeta\in \cl{\H}_{\C}$, we have (setting $z=\phi(\zeta$)):
     \begin{equation} \label{eqn:positionofG}
      |G(\zeta) - (\zeta+2)| =
            |\tilde{\phi}^{-1}(h(z)) - \tilde{\phi}^{-1}(f(z))|
              \leq \eps < 1, \end{equation}
     again using Lemma~\ref{lem:approximation} (\ref{item:approxinV}).

    In particular, $\re G^n(\zeta)\to \infty$ as $n\to\infty$; hence
     if $z\in V$, then the sequence $(g^n(z))$ lies in $V$ and has no accumulation point in
     $W$. By Lemma \ref{lem:accumulationsets}, there is a curve
     $\alpha\subset V$ whose accumulation set is exactly
     the accumulation set of the orbit of $z$. By choice of $V$,
     the accumulation set of $\alpha$ is exactly the accumulation set $K$
     of $\gamma$
     (recall Remark 2 after Lemma \ref{lem:strip}).

     Hence $V$ is contained in a Baker domain $U$ with
      $\omega(U)=K$, as required.
   \end{proof}

 \begin{proof}[Proof of Proposition \ref{prop:speed}]
  We now show how to modify the proof of Theorem \ref{thm:mainbaker}
   in order to obtain the extra claim in Proposition \ref{prop:speed}.
   We use the same notation as in the preceding proof.

  First we show that the claim holds for all $w$ in the Baker domain $U$
  such that
   $g^k(z)\in V$ for some $k\geq 0$,
   provided the
   function $\delta(t)$ was chosen sufficiently small, depending on the
   sequence $(\eps_n)$. Then we indicate how to modify the construction
   to ensure that $V$ is an \emph{absorbing set} for $U$;
   i.e., that every point of $U$ is mapped to $V$ under iteration.

  First note that we may assume without loss of generality that
   $(\eps_n)$ is a strictly
   decreasing sequence (otherwise, we instead consider the
   sequence $\bigl(\frac{1}{n}+\max_{k\geq n}\eps_k\bigr)$.

  Fix $w_0 := \gamma(0)$ as a base point. If $w\in V$, then by
   (\ref{eqn:positionofG}) the point $G^n(\phi^{-1}(w))$ is contained
   in the disk of radius $n\eps$ around $\phi^{-1}(w)+2n$. It follows that
     \begin{equation} \label{eqn:limitedgrowth}
         \dist_V(w_0,g^n(w)) = \dist_{\H}(\phi^{-1}(w_0),G^n(\phi^{-1}(w)))
             = O(\log n)\text{ \;as } n\to \infty.\end{equation}

 On the other hand, we can clearly let
   $\delta(t)$ tend to zero sufficiently
  rapidly to ensure that, for all sufficiently large $n$ and all
   $v\in V$ with $d(v,\partial W)\leq \eps_n$,
    \[ \dist_V(w_0,v) > n. \]

  Now let $w$ be a point such that
    $g^k(w)\in V$ for some $k\in\N$. By (\ref{eqn:limitedgrowth}),
     we have
     $\dist_V(w_0, g^{n+k}(w) ) < n$ for sufficiently large $n$.
      Hence
     \[ d(g^{n+k}(w),\partial W) > \eps_{n} > \eps_{n+k}, \]
    as desired.

 \smallskip

  It remains to show that we can ensure that $V$ is an absorbing set for $U$.
   To do so, we will modify the construction of $f$ by also requiring that there are many points
   near the boundary of $V$ that are not in the Baker domain $U$. This will ensure that points
   of $U\setminus V$ have very large hyperbolic distance from the ``central'' part of $V$, from
   which it follows that the orbits of all points of $U$ must enter $V$ eventually.

  More precisely, after picking $V$, $B$ and $C$, we also pick a
   discrete set $Z\subset W\setminus (\cl{V}\cup B\cup C)$ such that the hyperbolic
   metric of the set
   $W' := W\setminus Z$ satisfies
    \begin{equation} \label{eqn:hypdistanceWprime}
     \dist_{W'}(\phi(2n), W\setminus V) \to \infty \text{ \;as }n\to\infty.\end{equation}
  Also let $\Delta\subset W$ be a disk whose
   closure is disjoint from the sets $\cl{V}$, $B$, $C$ and $Z$, and let
   $\tilde{\Delta}\subset\C$ be another disk with
   $\pi(\tilde{\Delta})\subset \Delta$.

   Instead of
    letting the set $A$ from the construction of $h$
    consist only of $\cl{V}$, $B$ and $C$, we now set
     \[ A := \cl{V}\cup B \cup C \cup \cl{\Delta} \cup Z. \]
    We define $f$ as before on the first three sets, and such that
    $\cl{f(\Delta)}\subset\tilde{\Delta}$ and $f(Z)\subset \tilde{\Delta}$.

   If we choose $\eps$ sufficiently small in our application of Lemma
    \ref{lem:approximation}, then the approximating map $h$ also has the above properties.
    Thus $g=\pi\circ h$ satisfies $\cl{g(\Delta)}\subset\Delta$, so
    $\Delta$ is contained in the basin of an attracting fixed point of
    $g$. Since also $g(Z)\subset \Delta$,
    all points of $Z$ are also attracted to this fixed point, and in particular
    $U\subset W'$.

   Recall that, by (\ref{eqn:positionofG}), 
    the orbit of the point $w_1 := \phi(0)$ satisfies
     $|\phi^{-1}(g^n(w_1)) - 2n| \leq n$ for all $n\geq 0$. Hence, for
    $n\geq 1$,
     the hyperbolic distance $\dist_V(g^n(w_1),\phi(2n))$ is bounded by
     some constant $C$. Furthermore, if $w\in U$, then
     $\dist_U(g^n(w),g^n(w_1))\leq \dist_U(w,w_1)$ by Pick's theorem.
     So if $w\in U$, then
    \[ \dist_U(\phi(2n),g^n(w))\not\to\infty \text{ \;as } n\to\infty. \]
     Hence, by (\ref{eqn:hypdistanceWprime}),
     we have $g^n(w)\in V$ for sufficiently large $n$, as claimed.
 \end{proof}

 \section{Logarithmic tracts} \label{sec:logtracts}

  \begin{proof}[Proof of Theorem \ref{thm:mainlogtracts}]
   If $X=\Ch$, then let us assume that $a\neq w_0$. (The case
    $a=w_0$ is similar, but actually slightly simpler.
   We comment on this at the
    end of the proof.)

  Once again we pick $\pi:\C\to X$ to be a projection when $X$ is a torus, and otherwise
   a M\"obius transformation with $\pi(\infty)=a$; we also require that $\pi$ is
   chosen such that $\pi(0)=w_0$. Again we additionally pick discrete sets $B$ and $C$
   as in the proof of Theorem~\ref{thm:mainwandering}.

   Use Lemma \ref{lem:strip} to pick
    a simply-connected domain $V$ around $\gamma$ whose closure is
    disjoint from $B$ and $C$, and let
    $\phi:\H\to V$ be the corresponding Riemann map. Additionally, let
    $U'$ be the component of $\pi^{-1}(U)$ containing $0$; then
    $U'$ is a Jordan domain in $\C$. Also let $O'\subset\C$ be a
    Jordan domain
    containing $\cl{U'}$, but sufficiently close to $U'$ to ensure that
    $\pi$ is still injective on $\cl{O'}$. (That this is possible follows
    e.g.\ from the
    ``plane separation theorem''; see
    \cite[Chapter VI, Theorem (3.1)]{whyburnanalytictopology}.) Finally,
    set $O := \pi(O')$. Set $K:=\sup_{z\in O'}|z|$.

   We now let $f:V\to O'\setminus\{0\}$
    be a universal covering with $f(\partial V\cap W)=\partial O'$. Since $\phi^{-1}(\gamma)$ is bounded
    away from $\partial \H$ (see Lemma~\ref{lem:strip}(e)),
    we may also assume that $f$ is chosen such that
    $f(\gamma)\subset U'$.  As usual, we denote the extension
    to $\cl{V}_W$ by $f$ also. Define $f$ on $B$ and $C$ as in the proof
    of Theorems \ref{thm:mainwandering} and \ref{thm:mainbaker}.

   Using Lemma \ref{lem:arakelian}, we approximate $f$ by an analytic
    function $h:W\to\C$ such that $|h-f|$ is uniformly bounded on
    $B\cup C$, and such that
     \begin{equation}
       |h(z) - f(z) | < \eps\cdot |f(z)| \label{eqn:logarithmicapprox}
     \end{equation}
    for $z\in\cl{V}$, where $\eps<1/2$ is chosen small enough to ensure that
     $\eps < \dist_{\C}(\partial U', f(\gamma))/K$ and
     $\eps < \dist_{\C}(\partial O', U')/K$. 

    Then we have $h(\gamma)\subset U'$ and $h(z)\neq 0$ for all $z\in V$.
     Furthermore, the component $T$ of
     $h^{-1}(U')$ containing $\gamma$ does not intersect $\partial V$, and
     hence is completely contained in $V$. By
     the maximum principle, $T$ is simply-connected.

    We can thus
     form a logarithm $H := \log h : T \to \C$. Likewise, we can take
     a logarithm $F := \log f: V\to \C$.
     Note that $F:V\to F(V)$ is a conformal isomorphism
     by definition, and that $H(T)\subset F(V)$.

   The approximation condition
     (\ref{eqn:logarithmicapprox})
     implies that
      there is some
     $C>0$ such that
        \begin{equation}
            |H(z) - F(z)| < C  \label{eqn:bounded}
        \end{equation}
     for all $z\in T$.
     Consequently it follows that $H(z)\to\infty$ as
      $z\to\partial T\cap\partial W$ in $T$.
      We also have $H(z)\to \partial U'$ as
      $z\to \partial T\cap W\subset V$.
      Hence $H$ is a proper map, and therefore has a degree $d\geq 1$.

    We claim that $d=1$. It suffices to show that $H$ has no
     critical points in $T$. So let $c\in T$, and
     pick some simple closed curve
     $\alpha\subset H(T)$ that does not pass through any of the
     -- finitely many -- critical values of $H$ and
     such that $\alpha$ surrounds both
     $H(c)$ and a disk $D(z_0,C)$ of radius $C$ around some point $z_0\in H(T)$.
     This is possible because $H(T)$ contains a left half plane.

    The component $\beta$ of $H^{-1}(\alpha)$ that surrounds $c$ is
     a simple closed curve.
     By choice of $\alpha$, and by
     (\ref{eqn:bounded}), the curves $H\circ\beta$ and
     $F\circ\beta$ (where we fix some parametrization of $\beta$)
     have the same winding number $N\neq 0$ around
     $z_0$. Since $F$ is
     a conformal isomorphism, we must have $N= 1$, so
     $c$ is not a critical point of $H$, as required.

    So we have seen that $H$ is
     a conformal isomorphism, and hence
     $h$ maps $T$ as a universal covering map over $U'$,
     as claimed. Because $\re \phi^{-1}(\gamma(t))\to\infty$ as $t\to\infty$, we also
     have $h(\gamma(t))\to 0$ as $t\to\infty$. Hence the map
     $g := \pi\circ h$ has the stated properties.

 \medskip

 In the case where $X=\Ch$ and $a=w_0$, the proof proceeds essentially as
  above except that now
  $\pi(\infty)=a=w_0$. In this case $U'$ and $O'$ are
  neighborhoods of $\infty$. Also, we require the
  approximation to satisfy
     \[
       |h(z) - f(z) | < \eps
     \]
  instead of (\ref{eqn:logarithmicapprox}); thus, we can actually apply
  Theorem~\ref{thm:arakelian} directly instead of using Lemma~\ref{lem:arakelian}.
\end{proof}

\section{Hyperbolic surfaces}\label{sec:hyperbolic}
 \begin{proof}[Proof of Theorem~\ref{thm:mainhyperbolic}]
  The proofs of the three statements in the case where $X$ is hyperbolic
   are entirely analogous to those in the case of genus $\leq 1$, except
   that here we choose the map $\pi$ to be a universal covering map from
   a subset of the complex plane to $X$; e.g.\ $\pi:\D\to X$. 

  Then the constructions proceed as before. The
   composition $g:=\pi\circ h$ is again an Ahlfors islands map and has
   the desired properties, but the domain of definition $\wt{W}\subset W$ 
   of $g$ is almost certainly not 
   connected. 

  To obtain also the claim in Theorem \ref{thm:mainhyperbolic} that
   $g$ can be chosen with $W'$ connected, we can simply restrict
   $g$ to the component of $\wt{W}$ containing $\gamma$ in the case
   of Baker domains or logarithmic tracts
   (Theorems \ref{thm:mainbaker}
   and \ref{thm:mainlogtracts}).
   However, in the case of wandering domains (Theorem \ref{thm:mainwandering}),
   it is 
   possible that different components of the orbit of the wandering 
   domain belong to different components of $\wt{W}$. To prevent
   this from happening, we need to modify the construction from the proof
   of Theorem \ref{thm:mainwandering} slightly. 

  This is easy to do: Choosing the sets $B$ and $C$ as well as the domains
    $\Delta_k$ as before, we connect the disks $\Delta_k$ by arcs to
    form a ``tree'' $T\subset W\setminus(B\cup C)$, 
    which is a Weierstra{\ss} set 
    whose interior coincides precisely with 
    $\bigcup_{k=1}^{\infty} \Delta_k$. We can then
    carry out the approximation in a way that ensures that
    $h(T)\subset\D$ (recall that $\D$ is the domain of definition of
    $\pi$), and hence such that all domains $\Delta_k$ are contained in the
    same component of the domain of definition of $g$.
 \end{proof}

\section{Proof of Lemma \ref{lem:strip}}\label{sec:strip}

   \begin{proof}[Proof of Lemma~\ref{lem:strip}]
  We begin by noting that we will construct $V$ with the additional property
   that
   \begin{enumerate}
    \item[(d')] If $\alpha\subset V$ is any curve that tends to $\infty$ in
     $\hat{W}$, then $d(\alpha,\gamma(t))<\delta(t)$ for all sufficiently large $t$.
   \end{enumerate}
  Clearly we may assume without loss of generality that $\delta$ is bounded from above, say
   $\delta(t)\leq 1$, so (d'), together with (a), indeed implies part (d) of the 
   lemma. Furthermore, we may assume (modifying the function $\delta$ as necessary)
  that
  $d$ is a natural conformal metric on $W$ (i.e., the hyperbolic metric if $W$ is
  hyperbolic, or the Euclidean metric if 
  $W=\C$ or $W=\C/\Z$). Let us denote balls with respect to this metric by
  $B_d(z,r)$. 

 We may furthermore assume for simplicity 
  that $\delta$ is non-increasing, and that
  the closed ball $\overline{B_d(\gamma(t),\delta(t))}$ is homeomorphic
  to the closed unit disk (i.e., the injectivity radius of $W$ at $\gamma(t)$ is
  larger than $\delta(t)$). 

 We now define a function $\hat{\delta}(t)$ as follows. Let $t_0\in [0,\infty)$,
  and let us define
   \begin{align*}
      t_- &:= \min\bigl\{t\in [0,t_0]: \gamma\bigl((t,t_0]\bigr)\subset
               B_d(\gamma(t_0),\delta(t_0)/2) \}\quad\text{and} \\
      t_+ &:= \max\bigl\{t\in[t_0,\infty): \gamma\bigl([t_0,t)\bigr)\subset
               B_d(\gamma(t_0),\delta(t_0)/2) \}. \end{align*}
  We define
   \[ \hat{\delta}(t_0) := \frac{1}{2}\cdot 
           \inf_{t\notin [t_-,t_+]} d(\gamma(t),\gamma(t_0))>0. \] 
  Note that $\hat{\delta}(t_0)\leq \delta(t_0)/4$ because
   $d(\gamma(t_+),\gamma(t_0))=\delta(t_0)/2$. 

 Let us set
  \[ V_1 := \bigcup_{t\geq 0} B_d(\gamma(t),\hat{\delta}(t)) \]
  and study the set $\wt{V_1}$, 
   where $\wt{U}$
   denotes the union of a set $U$ and its compact complementary components. 
  The construction ensures that 
  \begin{enumerate}
   \item[(1)] $\wt{V_1}$ is simply-connected;
   \item[(2)] $\wt{V_1}\subset \bigcup_{t\geq 0} \{z\in W:d(z,\gamma(t))<\delta(t)\}$;
   \item[(3)] if $t_3\gg t_2 \gg t_1$, then $\gamma(t_1)$ and $\gamma(t_3)$ belong
     to different connected
     components of $\wt{V_1}\setminus B_d(\gamma(t_2),\delta(t_2))$. 
  \end{enumerate}
 Property (3) implies, in particular, that
  $\wt{V_1}$---and hence any subset $V$ of $\wt{V_1}$---satisfies
  (d').
 Now let $\phi:\H\to \wt{V_1}$ be a conformal isomorphism, which we can assume normalized such
  that $\phi^{-1}(\gamma)$ accumulates at $\infty$. Then (3) ensures---e.g.\
  by using the theory of prime ends---that $\phi(z)\to\infty$ in $\hat{W}$ iff
  $z\to\infty$ in $\H$, and in particular that $\phi^{-1}(\gamma)$ accumulates only at $\infty$. 
  Thus the domain $\wt{V_1}$ has all the properties required
  in the statement of the theorem, apart possibly from (b). However, we can choose
  a domain $H\subset\H$ such that $\partial H\cap\C$ is a simple curve contained in $\H$
  and such that $\phi^{-1}(\gamma)\subset H$. Then $V:=\phi(H)$ has the desired properties. 

\smallskip

 Now let us suppose that the curve $\gamma$ is $C^1$, and use a different construction
  to obtain a domain also satisfying 
  (e). First we extend $\gamma$ to be
 a $C^1$ curve on $(-1,\infty)$ and reparametrise if necessary to ensure that $\gamma$
 is `unit speed'. For $t>-1$, let $I(t)$ denote the line segment with centre at $\gamma(t)$,
 length $2\eps(t)$ and normal to the curve $\gamma$ at $\gamma(t)$. Here
 $\eps:(-1,\infty)\to (0,\infty)$ is a sufficiently rapidly decreasing convex function.
 The fact that $\gamma$
 is locally uniformly $C^1$ implies that, if $\eps$ is chosen sufficiently small, then 
 \[
 V=\bigcup_{-1<t<\infty}I(t)
 \]
 forms a simply-connected tubular neighbourhood of $\gamma|_{[0,\infty)}$ in $W$.  It
 is easy to check that $V$ satisfies parts~(a),~(b),~(c) and~(d) if we choose the
 function $\eps$ small enough. 

 To prove part~(e), let $z^+(t)$ and $z^-(t)$ denote the endpoints of $I(t)$, and let
 \[\partial^+ V=\bigcup_{-1<t<\infty}z^+(t)\quad \text{and} \quad \partial^- V=\bigcup_{-1<t<\infty}z^-(t).\]

 We claim that, provided the function $\eps$ was chosen sufficiently small, 
  there exists an absolute constant $c$, $0<c<1$, such that, for all
 $t\ge 0$,
 \begin{equation}\label{e5}
 \omega(\gamma(t), \partial^+ V,V)\le c,\quad\text{and}\quad \omega(\gamma(t), \partial^- V,V)\le c,
 \end{equation}
 where $\omega(z,E,V)$ denotes the harmonic measure in $V$ of the set $E\subset
 \partial V$ at the point $z\in V$. The inequalities~(\ref{e5}) imply that for large
 $t$, the point $\phi^{-1}(\gamma(t))$ lies in a set of the form $\{\zeta: \re \,\zeta\ge C|\zeta|\}$,
 where $C>0$, which proves part~(e).

 To prove~(\ref{e5}) note that for each $t>0$, the line segment $I(t)$ lies on the symmetric
 axis of a right-angled isosceles triangle, $\Delta(t)$ say, with vertex at $z^+(t)$,
 height $4\eps(t)$ and base outside the domain $V$. Hence, by the maximum principle,
 $\omega(\gamma(t), \partial^+ V,V)$ is dominated by the harmonic measure in $\Delta(t)$
 of the union of the two equal-length sides of $\Delta(t)$ at $\gamma(t)$, and this value
 is an absolute constant in $(0,1)$. This proves~(\ref{e5}).
\end{proof}
\begin{remark}
 Note that the proof shows that, provided $\gamma$ is $C^1$, the domain $V$ can be chosen
  such that 
  $\phi^{-1}(\gamma(t))\to\infty$ \emph{nontangentially}. It is not difficult to see that, 
  without a regularity assumption on $\gamma$, this cannot always be achieved. 

 For the conclusion of the lemma, we only require that $\phi^{-1}(\gamma)$ eventually
  enters every right half plane; i.e.\ that it tends to $\infty$ \emph{horocyclically}. 
  This can be achieved under much weaker conditions; e.g.\ it is sufficient to
  assume that there is a sequence $t_j$ with $t_j\to\infty$ such that
  $\gamma$ is differentiable in $t_j$ for all $j$. 
\end{remark}

\bibliographystyle{hamsplain}
\bibliography{/Latex/Biblio/biblio}

\providecommand{\href}[2]{#2}\def\polhk#1{\setbox0=\hbox{#1}{\ooalign{\hidewid%
th \lower1.5ex\hbox{`}\hidewidth\crcr\unhbox0}}}
  \def\polhk#1{\setbox0=\hbox{#1}{\ooalign{\hidewidth\lower1.5ex\hbox{`}\hidew%
idth\crcr\unhbox0}}} \input{cyracc.def} \def\j{{\u i}} \def\J{{\u I}}
  \newfont{\cyrit}{wncyi10 at 12pt}\def\cprime{$'$}
\providecommand{\bysame}{\leavevmode\hbox to3em{\hrulefill}\thinspace}
\begin{thebibliography}{10}

\bibitem[BKZ]{baranskikarpinskazdunik}
Krzysztof Bara{\'n}ski, Bogus{\l}awa Karpi{\'n}ska, and Anna Zdunik,
  \emph{Hyperbolic dimension of {J}ulia sets of meromorphic maps with
  logarithmic tracts}, Int. Math. Res. Not. IMRN (2008), Art. ID rnn141, 10,
  \mbox{\href{http://www.arXiv.org/abs/0711.2672v1}{arXiv:0711.2672v1}}.

\bibitem[B1]{waltermero}
Walter Bergweiler, \emph{Iteration of meromorphic functions}, Bull. Amer. Math.
  Soc. (N.S.) \textbf{29} (1993), no.~2, 151--188,
  \mbox{\href{http://www.arXiv.org/abs/math.DS/9310226}{arXiv:math.DS/9310226}%
}.

\bibitem[B2]{walterahlfors}
\bysame, \emph{A new proof of the {A}hlfors five islands theorem}, J. Anal.
  Math. \textbf{76} (1998), 337--347.

\bibitem[E1]{adamthesis}
Adam~L.~Epstein, \emph{Towers of finite type complex analytic maps}, Ph.D.~thesis, City
  Univ.~of New York, 1995, available at
  \url{http://pcwww.liv.ac.uk/~lrempe/adam/thesis.pdf}.


\bibitem[E2]{finitetype}
\bysame, \emph{Dynamics of finite type complex analytic maps {I}:
  Global structure theory}, Manuscript.

\bibitem[EO]{epsteinoudkerkahlfors}
Adam Epstein and Richard Oudkerk, \emph{Iteration of {A}hlfors and {P}icard functions which
  overflow their domains}, Manuscript.

\bibitem[G]{gaier}
Dieter Gaier, \emph{Vorlesungen \"uber {A}pproximation im {K}omplexen},
  Birkh\"auser Verlag, Basel, 1980.

\bibitem[M]{jackdynamicsthird}
John Milnor, \emph{Dynamics in one complex variable}, third ed., Annals of
  Mathematics Studies, vol. 160, Princeton University Press, Princeton, NJ,
  2006.

\bibitem[O]{oudkerkexoticbakerdomains}
Richard Oudkerk, \emph{{A}hlfors functions with exotic
  {B}aker domains}, Manuscript, 2004.

\bibitem[R]{hypdim}
Lasse Rempe, \emph{Hyperbolic dimension and radial {J}ulia sets of
  transcendental functions}, 2007,
  \mbox{\href{http://www.arXiv.org/abs/0712.4267}{arXiv:0712.4267}}.

\bibitem[RS]{ripponstallardslowescaping}
Philip~J.~Rippon and Gwyneth~M.~Stallard, \emph{Slow escaping points of
  meromorphic functions}, Preprint, 2008,
  \mbox{\href{http://www.arXiv.org/abs/0812.2410}{arXiv:0812.2410}}.

\bibitem[S]{scheinbergapproximation}
Stephen Scheinberg, \emph{Uniform approximation by functions analytic on a
  {R}iemann surface}, Ann. of Math. (2) \textbf{108} (1978), no.~2, 257--298.

\bibitem[W]{whyburnanalytictopology}
Gordon~T.~Whyburn, \emph{Analytic {T}opology}, American Mathematical
  Society Colloquium Publications, v. 28, American Mathematical Society, New
  York, 1942.

\end{thebibliography}

\end{document}